\documentclass[11pt,twoside,a4paper]{amsart}
\usepackage[inner=2cm,outer=2cm,top=2.5cm,headsep=1cm, footskip=1cm,bottom=3cm]{geometry}

\usepackage{amssymb, amscd, amsmath, enumerate,verbatim}
\usepackage{tikz}
\usepackage{color}
\usepackage{graphicx}
\usepackage[utf8]{inputenc}

\input xy
\xyoption{all}

\numberwithin{equation}{section}
\newtheorem{theorem}{Theorem}[section]
\newtheorem{lemma}[theorem]{Lemma}

\theoremstyle{definition}

\newtheorem{example}[theorem]{Example}

\newtheorem{definition}[theorem]{Definition}

\newcommand{\pb}{\ar@{}[dr]|{\mbox{\LARGE{$\lrcorner$}}}} %pull backs
\usepackage{caption}

\newcommand{\mk}{\mathbf{k}}

\newcommand{\D}{\mathcal{D}}

\newcommand{\dirlim}{\underset{\longrightarrow}{\mathrm{colim}}\,}

\newcommand{\id}{\mathrm{id}}

\newcommand{\Ass}{\mathcal Ass}

\newcommand{\Mag}{\mathcal Mag}

\begin{document}

\title{Up-to-homotopy algebras with strict units \newline (Extended abstract version)}

\author{Agust\'{\i} Roig}
\address[A. Roig]{Dept. Matemàtiques \\ Universitat Polit\`{e}cnica de Catalunya, UPC  \\ Diagonal 647, 08028 Barce\-lo\-na, Spain. https://orcid.org/0000-0002-4589-8075 }

\thanks{Roig partially supported by projects MTM2015-69135-P(MINECO/FEDER), 2014SGR634 and 2017SGR932}

\keywords{Up-to-homotopy algebra, operad, operad algebra, minimal model}
\subjclass[2010]{18D50, 55P62}

\begin{abstract} We prove the existence of minimal models à la Sullivan for operads with non-trivial arity zero. So up-to-homotopy algebras with strict units are just operad algebras over these minimal models. As an application, we give another proof of the formality of the \emph{unitary} $n$-little disks operad over the rationals.
\end{abstract}

\date{\today}
\maketitle
\tableofcontents

\section{Introduction and main result}

In the beginning, in Stasheff's seminal papers \cite{Sta63}, $A_\infty$-spaces (algebras) had points (units) in what was subsequently termed the zero arity of the unitary associative operad $\Ass_+$. They were still present in \cite{May72} and \cite{BoVo73}, for instance, but after that, points or units generally disappeared, and for a while, people working with operads assumed as a starting point $P(0) = \emptyset$, in the topological setting, or $P(0) = 0$ in the algebraic one: see for instance \cite{GK94}. This may have been caused because of the problems posed by those points (units); see \cite{Hin03}, or \cite{Bur18} for two examples, or, more to the point, \cite{Mar96} (as well as \cite{MSS02}), where Markl constructs minimal models for operads of chain complexes over a field of zero characteristic, carefully excluding operads with non-trivial arity zero.

More recently, the situation changed, and people have turned their efforts to problems involving non-trivial arity zero. In the topological context, for instance, we have \cite{MT14}; in the algebraic context, \cite{FOOO09a}, or \cite{HM12}; and dealing with both \cite{Fre17a}.

When introducing points (units) back in the theory of up-to-homotopy things, there are two main possibilities: either you consider \emph{strict} ones, as in Stasheff's original papers \cite{Sta63}, or in \cite{May72}, \cite{Fre17a}, \cite{Bur18}, or you consider \emph{up-to-homotopy} ones, or other relaxed versions of them: \cite{BoVo73}, \cite{HM12}, \cite{MT14}, and many others. You can even do both: \cite{KS09}.

In this paper, we work in the algebraic and strict part of the subject. The contribution we add to the present panorama is to prove the existence of Sullivan minimal models $P_\infty$ for operads $P$ on cochain complexes over a characteristic zero field $\mk$, with non-trivial arity zero in cohomology, $HP(0) = \mk$. In doing so, we extend the work of Markl \cite{Mar96}, (see also \cite{MSS02}) which proved the existence of such models for non-unitary operads, $P(0) = 0$. Our models include the one of \cite{Bur18} for the unitary associative operad $\Ass_+$. More precisely, our main result says:

\begin{theorem}
Every cohomologically unitary, $HP(0) = \mk$, cohomologically connected, $HP(1) = \mk$, and with a unitary multiplication operad $P$, has a Sullivan minimal model $P_\infty \stackrel{\sim}{\longrightarrow} P$.
\end{theorem}

In the non-unitary case, the importance of such minimal models is well known. For instance, they provide a \emph{strictification} of up-to-homotopy algebras, in that for an operad $P$ (with mild hypotheses), up-to-homotopy $P$-algebras are the same as strict, regular $P_\infty$-algebras. We show how up-to-homotopy associative algebras or $A_\infty$-algebras with strict units are exactly $(\Ass_+)_\infty$-algebras. As a second application, we offer another proof of the formality of the \emph{unitary} $n$-little disks operad $\D_{n+}$ over the rationals. This fills the gap in our paper \cite{GNPR05} noticed by Willwacher in his speech at the 2018 Rio International Congress of Mathematicians \cite{Wil18}.

\section{Ingredients}

Our result has been made possible thanks to two main ingredients: (1) the recently introduced $\Lambda$-modules and $\Lambda$-operads, of \cite{Fre17a}, and (2) the simplicial and Kan-like structures we found in an operad with unitary multiplication. Let us explain their role.

\subsection{Restriction operations} Sullivan minimal models are constructed by a cell-attaching inductive algorithm. In their original context of commutative dg algebras, the building blocks of this algorithm are called \emph{Hirsch extensions} \cite{GM13}, or \emph{KS-extensions} \cite{Hal83}. In the context of operads, they are called \emph{principal extensions} \cite{MSS02}. Their definition in the non-unitary case is the following.

\begin{definition}\label{OpdefKS} (See \cite{MSS02}.) Let $n\geq 2$ be an integer. Let $P = \Gamma (M)$ be free as a graded operad, where $M $ is a graded $\Sigma$-module, with $M(0)=M(1)=0$. An \emph{arity $n$ principal extension} of $P$ is the free graded operad

$$
P \sqcup_d \Gamma (E) :=\Gamma (M \oplus E) \ ,
$$

\noindent where $E$ is an arity-homogeneous $\Sigma_n$-module with zero differential and $d:E\longrightarrow ZP(n)^{+1}$  a map of $\Sigma_n$-modules of degree $+1$. The differential $\partial$ on $P\sqcup_d \Gamma (E)$ is built upon the differential of P, $d$, and the Leibniz rule.
\end{definition}

A Sullivan minimal model of an operad $P$ is a colimit of such principal extensions.

\begin{definition}\label{defSullmin} Given an operad $P$, a \emph{Sullivan minimal model} is a quasi-isomorphism $\rho : P_\infty \stackrel{\sim}{\longrightarrow} P$, which is built inductively on the arity of the operad through consecutive principal extensions
\end{definition}

$$
\xymatrix{
{P_2 = \Gamma (E(2))}\ar[r]\ar[rrrrd]^{\rho_2}  &{\dots}\ar[r] &{P_n = P_{n-1}\sqcup_{d_n} \Gamma (E(n))}\ar[r]\ar[rrd]^{\rho_n} &  {\dots}  \ar[r]  & {\dirlim_n P_n = P_\infty} \ar[d]^{\rho} \\
                   &  &  &  &  {P}
}
$$

\noindent in such a way that, for all $n$, $\rho_n : P_n \longrightarrow P$ is a quasi-isomorphism up to arity $n$.

This works perfectly fine in Markl's non-unitary case. The success of the Sullivan algorithm relies on the fact that, when restricted to operads which are cohomologically non-unitary $HP(0) = 0$ and cohomologically connected $HP(1) = \mk$, their minimal model is a free graded operad $P_\infty = \Gamma (M)$ over a $\Sigma$-module $M =  \bigoplus_{n=2}^{\infty} E(n) $ which is trivial in arities $0$ and $1$, $M(0) = M(1) = 0$. As a consequence, $P_\infty (n) = P_n(n)$: generators added in arities $ > n$ don't change what we have in lower ones.

The problem in introducing units $1 \in \mk = HP(0)$ of our cohomologically unitary operads as generators in the arity zero of their minimal model $P_\infty = \Gamma (M)$  would be that units give rise to \emph{restriction operations} which lower the arity:

$$
\delta_i = \_ \circ_i 1 : P(n) \longrightarrow P(n-1) \ , \quad \omega \mapsto \delta_i (\omega) = \omega\circ_i 1 \quad i = 1,\dots , n .
$$

So, in the presence of units, new generators $\omega \in E(n)$ added in arity $n\geq 2$ would produce \emph{new} elements in lower arities $\delta_i (\omega) \in P_n(n-1)$. Therefore, we would be introducing new generators in arity $n$ that would change lower arities; that is, in the previous steps of the induction process thus ruining it.

Nevertheless, we can  also assume that the generating module $M$ also has trivial arities $0$ and $1$ in our unitary case. This possibility has been recently made feasible thanks to Fresse's $\Lambda$-modules and $\Lambda$-operads, \cite{Fre17a}: to put it succinctly, we strip out of the operad all the structure carried by the elements of $P(0)$ and add it to the underlying category of $\Sigma$-modules, obtaining the category of $\Lambda$-modules. This way, we obtain a substitute for the general free operad functor with the bonus of getting our field $\mk$ in arity zero, with no need of any generators in the risky arities zero and one.

But we must keep track of those units somewhere if we want to build minimal models for cohomologically \emph{unitary} operads. This is how we do it: we add the restriction operations $\delta_i$ to the building blocks of our minimal model, \emph{without producing new generators in} $P_\infty$, with a just slight modification of the principal extensions.

\begin{definition}\label{OpdefKS+} Let $n\geq 2$ be an integer. Let $P$ be free as a unitary graded operad, $P= \Gamma (M)$, where $M $ is a graded $\Sigma$-module, with $M(0)=M(1)= 0$. A \emph{unitary arity $n$ principal extension} of $P$ is the free graded operad

$$
P \sqcup_d^\delta \Gamma (E) :=\Gamma (M \oplus E) \ ,
$$

\noindent where $E$ is an arity-homogeneous $\Sigma_n$-module with zero differential and:

\begin{enumerate}
  \item[(a)] $d:E\longrightarrow ZP(n)^{+1}$  is a map of $\Sigma_n$-modules of degree $+1$; the differential $\partial$ on $P\sqcup_d \Gamma (E)$ is built upon the differential of P, $d$ and the Leibniz rule.
  \item[(b)] $\delta_i : E \longrightarrow P(n-1), i = 1, \dots , n$ are morphisms of $\mk$-graded vector spaces, compatible with $d$ and the differential of $P$, in the sense that, for all $i=1, \dots , n$ we have commutative diagrams

      $$
      \xymatrix{
      {E} \ar[r]^{d} \ar[d]^-{\delta_i}   &  {ZP(n)^{+1}} \ar[d]^{\delta_i}  \\
      {M(n-1)}   \ar[r]^-{\partial}       &   {P(n-1)^{+1}} \ .
      }
      $$

      \noindent Which also have to be compatible with the $\Lambda$-structure of $P$, from arity $n-1$ downwards.
\end{enumerate}
\end{definition}

\subsection{A Kan-like structure} However, once we put \emph{unitary} principal extensions, with their extra restriction operations, in the Sullivan inductive algorithm, we have a new problem. To extend our \lq\lq partial" quasi-isomorphism $\rho_{n-1}$ to the next arity $\rho_n$, we now need to check that it is compatible with these new restriction operations.

For this, we introduce  \emph{simplicial-like} and \emph{Kan-like} structures in an operad with unitary multiplication. To the best of our knowledge, both structures are new.

We begin with the simplicial-like structure. Restriction operations $\delta_i$ give us the \emph{face maps}. To obtain the \emph{degeneracy maps}, we need a unitary multiplication on $P$; that is, a morphism of operads $\Mag_+ \longrightarrow P$. Here, $\Mag_+$ stands for the \emph{unitary magmatic operad}. $\Mag_+$ is the operad of unitary magmas: algebras with a unit, a single operation, and just the unit relations. This morphism gives us elements $1\in P(0)$ and $m_2 \in P(2)$. These elements behave as a unit and a product: $m_2 \circ_1 1 = \id = m_2 \circ_2 1$. This extra condition of a unitary multiplication has an easy and natural interpretation: we are only asking that the unit of our operad $1 \in P(0)$ not be an \lq\lq idle" one: there needs to be an arity two operation $m_2$ for which $1$ \emph{actually} works as a unit.

With this unitary multiplication in $P$, we define our degeneracy maps as

$$
s_i : P(n) \longrightarrow P(n+1) \ , \quad s_i(\omega) = \omega\circ_i m_2 \ , \quad i = 1, \dots , n \ .
$$

It's an easy exercise to check that these $\delta_i$ and $s_i$ fulfill the simplicial identities needed to prove lemma \ref{Kan-likeresult} below. Moreover, since each $P(n)$ is an abelian group, $P$ would look like a sort of \emph{Kan complex} (up to a shift). Nevertheless, this is \emph{not} the Kan structure we are interested in, but the following:

\begin{definition}\label{Kan-likecondition} Let $\left\{\omega_i \right\}_{i=1, \dots , n}$ be a family of elements in $P(n-1)$. We say that they verify a \emph{Kan-like condition} if $\delta_i \omega_j = \delta_{j-1}\omega_i$, for all $i < j $.
\end{definition}

\begin{example} Elements $\omega \in P(n), \ n \geq 1 ,$ produce families $\left\{ \omega_i = \delta_i \omega\right\}_{i = 1, \dots , n}$ in $P(n-1)$ that verify the \emph{Kan-like condition}.
\end{example}

The reciprocal of this example is also true.

\begin{lemma}\label{Kan-likeresult} Let $\left\{ \omega_i  \right\}_{i = 1, \dots , n}$ be a family of elements in $P(n-1)$ verifying the Kan-like condition. Then there exists an $\omega \in P(n)$ such that $\delta_i \omega = \omega_i$ for all $i = 1, \dots , n$.
\end{lemma}

Moreover, we can prove that, if all the $\omega_i$ are cocycles, coboundaries, or belong to the kernel or the image of an operad morphism $\varphi : P \longrightarrow Q$, then $\omega$ can be chosen to be also a cocycle, a coboundary, or to belong to the kernel or the image of $\varphi$, respectively. Even more: if the $\omega_i = \omega_i (e)$ depend linearly on $e \in E(n)$, we can choose $\omega = \omega (e)$ to depend $\Sigma_n$-equivariantly on $e$.

So much for the explanations. Let us now see all these constructions actually in action.

\begin{proof}[Sketch of the proof of theorem 1] To build $\rho_2 : P_2 \longrightarrow P$ in the non-unitary case, we take the generators in arity two to be $E=E(2) = HP(2)$. Then, we choose a $\Sigma_2$-equivariant section $s_2 : HP(2) \longrightarrow ZP(2) \subset P(2)$ of the projection $\pi_2 : ZP(2) \longrightarrow HP(2)$. And we get our first stage of the inductive algorithm as:

$$
P_2 = \Gamma (E) \ , \qquad \partial_{2 | E} = 0 \ , \qquad \text{and} \qquad \rho_2 : P_2 \longrightarrow P \ , \qquad \rho_{2 | E} = s_2 \ .
$$

In the unitary case, our section should make the following diagram to commute too:

$$
\xymatrix{
{E =HP(2)} \ar[d]^{\delta_i} \ar@{.>}@/^1pc/[r]^{s'_2}   &  {ZP(2)}  \ar[d]^{\delta_i} \ar[l]^-{\pi_2}     \\
{ \mk = HP(1)}  \ar@/^1pc/[r]^{s_1}    &  {ZP(1)} \ar[l]^-{\pi_1}
}
$$

Here, section $s_1$ is the unique $\mk$-linear map sending $\id \in HP(1)$ to $\id \in ZP(1)$ and the restriction operations on $E$ are the ones induced by $\delta_i : P(2) \longrightarrow P(1), i = 1,2$ on cohomology.

This is not necessarily true for the section $s_2$ we have found in the non-unitary case. So, given $e\in E$, we study the differences $\omega_i (e) = \delta_i s_2(e) - s_1\delta_i (e) \in P(1),  i=1,2$. And we observe that they are coboundaries and verify our  \emph{Kan-like} condition in definition \ref{Kan-likecondition}. Therefore, because of lemma \ref{Kan-likeresult}, we get a coboundary $\partial \omega (e) \in P(2)$, such that $\delta_i \partial \omega (e) = \omega_i(e), i=1,2$. With this, we modify the section $s_2$ from the non-unitary case to a new one $s'_2 (e) = s_2(e) - \partial\omega (e)$ which is compatible with the restriction operations $\delta_i$. Finally, we average over $\Sigma_2$

\begin{align*}
\widetilde{s}_2(e) &= \frac{1}{2!}\sum_{\sigma \in \Sigma_2} \sigma \cdot s'_2(\sigma^{-1}\cdot e)  \\
                   &= \frac{1}{2} \left( s'_2(e) + (2\ 1)\cdot s'_2((2\ 1)\cdot e) \right)
\end{align*}

\noindent and obtain a $\Sigma_2$-equivariant section, without losing anything we previously had for $s_2$ and $s'_2$. Therefore, we have our induced morphism of \emph{unitary} operads $\rho_2 : P_2 \longrightarrow P$.
\end{proof}

%    Numbered section headings
%\section{}
%\subsection{}

%    Unnumbered section headings
%\section*{}
%\subsection*{}

%    Ordinary theorem and proof
%\begin{theorem}[Optional addition to theorem head]
% text of theorem
%\end{theorem}

%\begin{proof}[Optional replacement proof heading]
% text of proof
%\end{proof}

%    Figure insertion
%\begin{figure}
%\includegraphics{filename}
%\caption{text of caption}
%\label{}
%\end{figure}
% You can also use a package such as \epsfig or \psfig.
% Please use \epsfig rather than \epsf, which has
% become obsolete.

%    Mathematical displays

% Numbered equation
%\begin{equation}
%\end{equation}

% Unnumbered equation
%\begin{equation*}
%\end{equation*}

% Aligned equations
%\begin{align}
% &  \\
% &
%\end{align}

%    Bibliography.


\begin{thebibliography}{99}

%\bibitem{}
\bibitem[BoVo73]{BoVo73} J. M. Boardman, R. M. Vogt, \emph{Homotopy invariant algebraic structures on topological spaces}, Lecture Notes in Mathematics, Vol. 347, Springer-Verlag, Berlin-New York (1973).
\bibitem[Bur18]{Bur18} J. Burke, \emph{Strictly unital {$\rm A_\infty$}-algebras}, J. Pure Appl. Algebra 222, (2018), 4099--4125.
\bibitem[FOOO09a]{FOOO09a} K. Fukaya, Y.-G. Oh, H. Ohta, K. Ono, \emph{Lagrangian intersection {F}loer theory: anomaly and obstruction. {P}art {I}}, AMS/IP Studies in Advanced Mathematics 46, American Mathematical Society, Providence, RI; International Press, Somerville, MA (2009).
\bibitem[Fre17a]{Fre17a} B. Fresse, \emph{Homotopy of operads and {G}rothendieck-{T}eichmüller groups. {P}art 1.
     {T}he algebraic theory and its topological background}, Mathematical Surveys and Monographs 217, American Mathematical Society, Providence, RI (2017).
\bibitem[GK94]{GK94} V. Ginzburg, M. Kapranov, \emph{Koszul duality for operads}, Duke Math. J. 76 (1994), 203--272.
\bibitem[GM13]{GM13} P. Griffiths, J. Morgan, \emph{Rational homotopy theory and differential forms}, second edition, Springer, New York (2013).
\bibitem[GNPR05]{GNPR05} F. Guillén, V. Navarro, P. Pascual, A. Roig, \emph{Moduli spaces and formal operads}, Duke Math. J. 129 (2005), 291--335.
\bibitem[Hal83]{Hal83} S. Halperin, \emph{Lectures on minimal models}, M\'em. Soc. Math. France (N.S.) 9-10 (1983).
\bibitem[Hin03]{Hin03} V. Hinich, \emph{rratum to \lq\lq{H}omological algebra of homotopy algebras"}, preprint arXiv:0309453 (2003).
\bibitem[HM12]{HM12} J. Hirsh, J. Mill\`es, \emph{Curved {K}oszul duality theory}, Mathematische Annalen 354 (2012), 1465--1520.
\bibitem[KS09]{KS09} M. Kontsevich, Y. Soibelman, \emph{Notes on {$A_\infty$}-algebras, {$A_\infty$}-categories and non-commutative geometry} in \emph{Homological mirror symmetry}, Lecture Notes in Phys. 757 (2009), 153--219.
\bibitem[May72]{May72}, J. P. May, \emph{The geometry of iterated loop spaces}, Lectures Notes in Mathematics, Vol. 271, Springer-Verlag, Berlin-New York (1972).
\bibitem[Mar96]{Mar96} M. Markl, \emph{Models for operads}, Comm. Algebra 24 (1996), 1471--1500.
\bibitem[MSS02]{MSS02} M. Markl, S, Shnider, J. Stasheff, \emph{Operads in algebra, topology and physics}, Mathematical Surveys and Monographs, American Mathematical Society, Providence, RI (2002).
\bibitem[MT14]{MT14} F. Muro, A. Tonks, \emph{Unital associahedra}, Forum Math. 26 (2014), 593--620.
\bibitem[Sta63]{Sta63} J. Stasheff, \emph{Homotopy associativity of {$H$}-spaces. {I}, {II}}, Trans. Amer. Math. Soc. 108 (1963), 275-292; ibid. (1963), 293--312.
\bibitem[Wil18]{Wil18} T. Willwacher, \emph{Little disks operads and {F}eynman diagrams} in \emph{Proc. Int. Cong. Of Math.} (2018).
\end{thebibliography}
\end{document}